\documentclass[11pt]{article}

\usepackage{amsmath,amsfonts,amssymb} 
\usepackage[english]{babel}                        
\usepackage{geometry}
\usepackage[T1]{fontenc}			                                

\usepackage[english]{babel}                        
\usepackage[T1]{fontenc}			   

\usepackage{theorem}{\theorembodyfont{\slshape}
  \newtheorem{theorem}{Theorem}
  \newtheorem{proposition}{Proposition}
  \newtheorem{corollary}{Corollary} 
  \newtheorem{lemma}{Lemma}
}{\theorembodyfont{\upshape}
  \newtheorem{remark}{Remark} 
}
\def\Proof{{\smallskip\noindent{\em Proof. }}}     
\def\endProof{{\hfill$\Box$\medskip\noindent}}     

\def\endProof{{\hfill$\Box$}}
\newcommand\R{{\mathbb{R}}}
\renewcommand\P{{\mathbb{P}}}
\newcommand\N{{\mathbb{N}}}

\newcommand\K{{\mathbb{K}}}

\newcommand\KK{{\mathfrak{K}}}

\renewcommand\S{{\mathbb{S}}}
\renewcommand\div{{\rm div}}

\renewcommand\L{{\mathbb{L}}}

\title{On the effect of external forces
on incompressible fluid motions at large distances}

 \author{
 Hyeong-Ohk Bae\\
 Departement of Mathematics, School of Natural Sciences\\
 Ajou University, Suwon 443--749, Korea\\
 {\tt hobae@ajou.ac.kr}\\
\\
  Lorenzo Brandolese\\
 Universit\'e de Lyon; Universit\'e Lyon 1\\
 CNRS, UMR 5208 Institut Camille Jordan,\\
 43, bd. du 11 novembre 1918\\
 F - 69622 Villeurbanne Cedex, France.\\
 {\tt brandolese@math.univ-lyon1.fr}\\
 }

\date{April 1st, 2009}
\begin{document}

\maketitle



\section
{Introduction}
\label{sec1}
Let $d\ge2$ be an integer.
We consider the Navier--Stokes equations in $\R^d$,
for an incompressible fluid submitted to an external force 
$f\colon\R^d\times \R^+\to\R^d $:
\begin{equation}
\label{nsd}
\begin{cases}
 \partial_t u+u\cdot \nabla u=\Delta u -\nabla p+f& \\
 \div{u}=0 &\\
 u(x,0)=a(x)
 \end{cases}
\end{equation}

We want to investigate the effect of the external force
on the large time and spatial asymptotics of the solution.
 In particular, we will show that even if $f$ is small and
well localized (say, compactly supported in space-time) but with
non-zero mean,
then the velocity of the fluid particles at \emph{all times} $t>0$
and in \emph{all point} 
$x$ outside of balls $B(0,R(t))$ with large radii,
is consirerably faster than in the case of the free Navier--Stokes
equations.

This sharply contrasts with the asymptotic properties of the
solution of other semilinear parabolic equations, where localized external forces do not affect
the behavior of the solution as $|x|\to\infty$:
only the large time behavior is influenced by the force.

The main results  of the paper will be stated in Section~\ref{sectionmain}:
they include new  sharp pointwise estimates
of the form
\begin{equation}  
 \label{uuop}
ct|x|^{-d}\le |u(x,t)|\le c't|x|^{-d},
\end{equation}
valid for $0<t<t_0$ small enough and $|x|\ge R(t)$.
Such pointwise estimates should be compared to those in the case of the  Navier--Stokes equations
without forcing studied in~\cite{BraV07}, where the behavior of $u$ is like $ t|x|^{-d-1}$.
The constant~$c$ in~\eqref{uuop} essentially depends on the integral of~$f$.

Estimates~\eqref{uuop} can be applied, {\it e.g.\/} to the case of compactly supported initial data.
In such case, they  describe how fast fluid particles start their motion in the far field,
at the beginning of the evolution. An even more precise description of the motion
of flows with localized data at large distances
will be given by the following asymptotic profile, valid for all $t>0$ such that the strong solution~$u$
is defined and for all $|x|\ge R(t)$:
\begin{equation}
 \label{befo}
u(x,t)\simeq \mathfrak{K}(x)\int_0^t\!\!\!\int f(y,s)\,dy\,ds,
\end{equation}
where $\mathfrak{K}(x)$ is the matrix of the second order derivatives of the fundamental solution
of the Laplacian in~$\R^d$. As such, $\mathfrak{K}_{j,k}(x)$ is a homogeneneous function of degree~$-d$.
See Theorem~\ref{theo1}, Section~\ref{sectionmain}, for a more precise statement.

As a consequence of our asymptotic profiles one can also recover some known bounds
on the large time behavior of $L^p$-norms for strong solutions
with  very short proofs.
See, {\it e.g.\/}, \cite{AGSS, BaeJ, choe-jin, HeX01, KuT07, OlivTiti, Schon85, Sch91, Skala09, Zhang04, Zhou07}
for a small sample of recent works on this  topic and related developements.
Our approach is different from that of these papers since it consists in deducing information on spatial norms
about the solution from information on their pointwise behavior.
As such, it is not so well suited for the study of  weak solutions as it is for that of strong, small solutions.
However, its advantage is that it allows to get sharp estimates, both
from above and below, also in the case of $L^p$-norms with weight like, {\it e.g.\/}, $(1+|x|)^\alpha$,
 for a wider range of the parameters~$p$ and~$\alpha$. Namely,
we will prove that, under localization assumptions on the datum and the force,
and provided 
$$\int_0^\infty\!\!\!\int f(y,s)\,dy\,ds\not=0_{\R^d},$$
then, for $1<p\le\infty$, $\alpha\ge 0$ and $t\to\infty$,
\begin{equation*}
\|(1+|x|)^\alpha u(t)\|_{p}\simeq  ct^{-\frac{1}{2}(d-\alpha-d/p)}, \qquad\hbox{if $ \alpha+d/p<d$}.    
\end{equation*}
The above restriction on the parameters~$p$ and $\alpha$ are optimal, as we will show that, for all $t>0$, and $1\le p<\infty$,
\begin{equation*}
 \|(1+|x|)^\alpha u(t)\|_{p}=+\infty \qquad\hbox{if $\alpha+d/p\ge d$}.
\end{equation*}
See Theorem~\ref{theo2}, Section~\ref{sectionmain}, for more precise statements.

\section{Preliminary material}

The assumptions on the external force will be the following:
\begin{eqnarray}
\label{f1}
&&|f(x,t)| \le \epsilon\Bigl[(1+|x|)^{-d-2}\wedge(1+t)^{-(d+2)/2}\Bigr]\\
\label{f2}
&&\|f\|_{L^1(\R^d\times\R^+)}\le \epsilon
\end{eqnarray}
for some (small) $\epsilon>0$ and a.e. $x\in\R^d$ and $t\ge0$.
Note that
assumption~\eqref{f2} can be interpretated
as a logarithmic improvement on the decay estimate~\eqref{f1}.
Of course, an additional potential force $\nabla\Phi$ could
also be added, affecting in this way the pressure of the fluid,
but not its velocity field.
On the other hand, because of condition~\eqref{f2}, and due to the unboundedness of singular integrals
(the Riesz transforms) in $L^1$ and, we are not allowed to apply the Helmoltz decompsition to~$f$.
Therefore, we will {\it not\/} restrict ourselves to divergence-free external forces, as it is sometimes the case
in the literature. 

The above assumptions on~$f$, as well as the smallness assumption on the datum below, look stringent:
we simply put assumptions that allow us to provide the  shortest possible
and self-contained construction of strong solutions with some decay.
Indeed, our main goal
will be to show that even in the case of very nicely behaved
external forces (with non-zero mean) and data,  the solution will be badly behaved at infinity.
More precisely, {\it upper bounds\/} on~$f$ and $a=u(0)$, no matter how good, will lead to
{\it lower bounds\/} on~$u(t)$.

But for the time being, we concentrate on the simpler problem of upper bounds on~$u$,
and start by establishing the following simple result.

\begin{proposition}
\label{prop1}
Assume that $f$ satisfies~\eqref{f1}-\eqref{f2}.
Let also $a\in L^1(\R^d)$ be a divergence-free vector field such that,  $\|a\|_1<\epsilon$
and $|a(x)|\le \epsilon(1+|x|)^{-d}$
for a.e. $x\in\R^d$.
If $\epsilon>0$ is small enough, then there exists a unique strong solution $u$
of the Navier-Stokes equation (NS)
satisfying, for some constant $C>0$, the pointwise decay estimates
\begin{equation}
\label{decu}
|u(x,t)|\le C\epsilon\Bigl[(1+|x|)^{-d}\wedge(1+t)^{-d/2}\Bigr],
\end{equation}
and such that $u(0)=a$ (in the sense of the a.e. and distributional convergence as $t\to0$).
In particular, for $\alpha\ge0$, $1<p\le \infty$ , we have
\begin{equation}
\label{ubu}
\Bigl\|(1+|x|)^{\alpha}u(t)\Bigr\|_{p}\le C(1+t)^{-\frac{1}{2}(d-\alpha-\frac{d}{p})},
\qquad
\hbox{when $\alpha+d/p<d$}.
\end{equation}
The above estimate remains true in the limit case $(\alpha,p)=(d,\infty)$.
\end{proposition}

\begin{remark}
The Proposition above can be viewed as the limit case of a
previous result by Takahashi \cite{Tak},
where pointwise decay estimates of the form 
$|u(x,t)|\le C\bigl[(1+|x|)^{-\gamma}\wedge (1+t)^{-\gamma/2}]$
had been obtained for $0\le \gamma<d$, 
with a different method and under different assumptions.

There are, on the other hand, many other methods for proving
sharp upper bounds of the form
$$\|(1+|x|)^\alpha u(t)\|_p\le Ct^{-\frac{1}{2}(d-\alpha-d/p)},.$$
For example, it would be possible to adapt the arguments of~\cite{choe-jin}
or~\cite{KuT07} to the case of non-zero external forces.
See also \cite{AGSS}, \cite{OlivTiti} \cite{Zhang04} for other
different approaches for getting decay estimates.
However, the optimal range of the parameters $\alpha$ and $p$ for the validity of such estimate
is not discussed in the previous papers.
Our condition $\alpha+d/p<d$ is more general, and turns out to be optimal whenever $f$ has non-zero mean,
as we will show in Theorem~\ref{theo2}.
\end{remark}

We denote by $\P$ be the usual Leray projector onto the solenoidal
vector fields and with $e^{t\Delta}$ the heat semigroup.
The solution of Proposition~\ref{prop1} is obtained by a straightforward fixed point argument by
solving the equivalent integral formulation:
\begin{equation}
\label{6a5}
\left\{
\begin{aligned}
&u(t)=e^{t\Delta}a-\int_0^t e^{(t-s)\Delta}\P\nabla\cdot(u\otimes u)(s)\,ds+
\int_0^t e^{(t-s)\Delta}\P f(s)\,ds &\\
&\div a=0.
\end{aligned}
\right.
\end{equation}

Let $\K(x,t)$ be the kernel of $e^{t\Delta}\P$.
It is well known, and easy to check that, for all $t>0$, $\K(\cdot,t)$ belongs
to $C^\infty(\R^d)$ and satisfies the scaling relation
$\K(x,t)=t^{-d/2}\K(x/\sqrt t,1)$, and the decay estimates
\begin{eqnarray}
\label{decK}
&&|\K(x,t)|\le C\Bigr[ |x|^{-d}\wedge t^{-d/2}\Bigr],\\
\label{denK}
&&|\nabla\K(x,t)|\le C\Bigr[ |x|^{-d-1}\wedge t^{-(d+1)/2}\Bigr].
\end{eqnarray}

The following Lemma is an immediate consequence of the above
decay estimates.

\begin{lemma}
\label{lemlog}
For all $(x,t)$ such that $|x|\ge e\sqrt t$ and $t>0$, the kernel $\K(x,t)$ 
satisfies:
\begin{equation}
\label{lnK}
\int_0^t\!\!\int_{|y|\le |x|}|\K(y,s)|\,dy\,ds\le Ct\log(|x|/\sqrt t).
\end{equation}
\end{lemma}

The proof of Proposition~\ref{prop1} relies
on the following lemma.

\begin{lemma}
\label{lemLf}
Let $f$ be such that \eqref{f1}-\eqref{f2} holds.
Let $\L$ be the linear operator defined by
\begin{equation}
\label{Lf}
\L(f)=\int_0^t e^{(t-s)\Delta}\P f(s)\,ds.
\end{equation}
Then there is a constant $C>0$ such that
\begin{equation}
\label{eLf}
|\L(f)(x,t)|\le C\epsilon\bigl[(1+|x|)^{-d}\wedge(1+t)^{-d/2}\Bigr].
\end{equation}
\end{lemma}

\Proof
We give a proof based on Lorentz spaces $L^{p,q}$, see~\cite{BL, Lem02} for their definition and basic
properties,
since Lorentz spaces will play a role also in the sequel.
For simplicity, we will drop the constant $\epsilon$ throughout the proof.
By inequality~\eqref{f1} we readily get
\begin{equation*}
\|f(t)\|_{L^{p,\infty}}\le C(1+t)^{-\frac{1}{2}(d+2-\frac{d}{p})}, \qquad 1<p<\infty.
\end{equation*}
Now, interpolating the $L^{p,q}$-space between $L^{p_1,\infty}$ and $L^\infty$,
with $1<p_1<p$ and $1\le q\le \infty$, we deduce
\begin{equation}
\label{depq}
\|f(t)\|_{L^{p,q}}\le C(1+t)^{-\frac{1}{2}(d+2-\frac{d}{p})}, \qquad 1<p<\infty,\quad 1\le q\le\infty.
\end{equation}
In particular, choosing $p=q=d$,
$$\|f(s)\|_{L^{d}}\le C(1+s)^{-(d+1)/2}.$$
On the other hand, it follows from~\eqref{decK} (or by well known $L^p$-estimates on the kernel~$\K$), that
$\|\K(\cdot,t-s)\|_{d/(d-1)}\le C(t-s)^{-1/2}$.
By Young inequality we get
$\|K(t-s)*f(s)\|_\infty \le C(t-s)^{-1/2}(1+s)^{-(d+1)/2}$.
Another obvious estimate estimate is
$$\|K(t-s)*f(s)\|_\infty \le C(t-s)^{-d/2}\|f(s)\|_1.$$
Now splitting the intergral defining~$\L$ at $s=t/2$ we get
\begin{equation*}
\|L(f)(t)\|_\infty \le C(1+t)^{-d/2}.
\end{equation*}
In particular, for $|x|\le e\sqrt t$, we obtain
$|\L(f)(x,t)|\le C(1+|x|)^{-d}$
and now we can limit ourselves to the case $|x|\ge e\sqrt t$ and $|x|\ge1$.
We  write $\L(f)=\L_1(f)+\L_2(f)$, where
$\L_1(f)=\int_0^t\!\!\int_{|x-y|\le |x|/2}\dots$ and
$\L_2(f)=\int_0^t\!\!\int_{|x-y|\ge |x|/2}\dots$. 
For treating $\L_1(f)$, we use assumption~\eqref{f1}
and obtain, after changing the variables $z=x-y$ and $\tau=t-s$,
by Lemma~\ref{lemlog},
\begin{equation*}
|\L_1(f)(x,t)|  \le Ct|x|^{-d-2}\log(|x|/\sqrt t)\le C|x|^{-d},
\end{equation*}
where  we used that $|x|\ge e\sqrt t$.
Using  estimate~\eqref{decK} and assumption~\eqref{f2}
we get $|\L_2(f)(x,t)|\le C|x|^{-d}$.

\endProof 

\bigskip
Direct estimates on the heat kernel show that for $a\in L^1(\R^n)$ satisfying
$\|a\|_1\le\epsilon$ and
$|a(x)|\le \epsilon(1+|x|)^{-d}$, we have,  for some constant $C>0$ independent on $x$ and $t\ge0$,
\begin{equation}
\label{rihs} 
|e^{t\Delta}a(x)|\le C\epsilon\Bigl[(1+|x|)^{-d}\wedge(1+t)^{-d/2}\Bigr].
\end{equation}
Owing to estimate~\eqref{eLf}, 
the fixed point argument argument used in Miyakawa~\cite{Miy00},
where the solution was constructed
in the special case $f\equiv0$,
goes through in our situation.
More precisely, the approximate solutions $u_k$,
constructed in the usual manner (see, {\it e.g.\/} \cite{Ka84, Miy00})
converge in the Banach space of measurable functions of the $(x,t)$-variables
bounded by the right-hand side of~\eqref{rihs}.
The existence and the unicity of a solution~$u$ satisfying
the pointwise estimates~\eqref{decu} follows.

It is straightforward to see that any function satisfying estimates~\eqref{decu}
must also verify the weighted-$L^p$ bounds~\eqref{ubu}.
This can be proved first by bounding the $L^{p,\infty}$-norms of $(1+|x|)^\alpha u(t)$ and then
interpolating the $L^p$-space between $L^{p_1,\infty}$ and $L^\infty$
with $1<p_1<p$. This is exactly the same argument that we applied to deduce estimate~\eqref{depq}.

This completes the proof of Proposition~\ref{prop1}.

\endProof

\bigskip

Denote by $\Gamma$ the Euler Gamma function and by $E_d$ the fundamental
solution of $-\Delta$ in $\R^d$.
The following Lemma (proved in~\cite{bra-fine})
will be useful:
\begin{lemma}
\label{asG}
Let $\KK=(\KK_{j,k})$, where $\KK_{j,k}(x)$ is the homogeneous function of degree~$-d$
\begin{equation}
\label{KK}
\KK_{j,k}(x)=
\partial^2_{x_j,x_k}E_d(x)=
\frac{\Gamma(d/2)}{2\pi^{d/2}}\cdot \frac{\bigl(-\delta_{j,k}|x|^2+dx_jx_k\bigr)} {|x|^{d+2}}.
\end{equation}
Then the following decomposition holds:
\begin{equation}
\label{profKF}
\K(x,t)=\KK(x)+|x|^{-d}\Psi\Bigl(x/\sqrt t\Bigr),
\end{equation}
where $\Psi$ 
 is a function defined on $\R^d$,
smooth outside the origin and such that, for all $\alpha\in\N^d$, and $x\not=0$, 
$|\partial^\alpha \Psi(x)|\le Ce^{-c|x|^2}$, 
where $C$ and $c$ are positive constant,
depending on $|\alpha|$ but not on $x$.

\end{lemma}

Before going further, let us put an additional assumption on the external force,
namely
\begin{equation}
\label{f3}
\|\,|x|f(\cdot,t)\|_{1}\le C(1+t)^{-1/2}, \qquad \hbox{for some $C>0$ and all 
$x\in\R^d$, $t\ge0$}.
\end{equation}
Just like condition~\eqref{f2}, we can view assumption~\eqref{f3}
as another logarithmic improvement on the decay estimate~\eqref{f1}.
However, here the smallness of the constant $C>0$ is unessential.

Next Lemma provides an explicit asymptotic expansion 
for $\L(f)(x,t)$, as $|x|\to\infty$. The function of the $(x,t)$-variables
$\mathcal{R}_f$ defined by relation~\eqref{asLf} below can be viewed
as a remainder term, that can be neglected at large distances, when
$|x|>\!\!\!>\sqrt t$.

\begin{lemma}
Let $f$ such that assumptions~\eqref{f1}-\eqref{f2} and \eqref{f3} hold.
For $x\not=0$, and $t>0$, define $\mathcal{R}_f(x,t)$, through the relation
\begin{equation}
\label{asLf}
\L(f)(x,t)=\KK(x) \int_0^t\!\!\int f(y,s)\,dy\,ds+\mathcal{R}_f(x,t).
\end{equation}
Then  $\mathcal{R}_f(x,t)$ satisfies, for some constant~$C>0$,
\begin{equation}
\label{br1}
\forall\, |x|\ge e\sqrt t\colon \quad
|\mathcal{R}_f(x,t)|\le C|x|^{-d-1}\sqrt{t}.
\end{equation}
\end{lemma}

\begin{remark}
\label{rema2}
The product on the right hand side of~\eqref{asLf} is the
usual product between the $d\times d$ matrix $\KK(x)$
and the $d$-vector field $\int_0^t\!\!\int f_k(y,s)\,dy\,ds$.
The $j$-component is thus given by $\sum_{k=1}^d \KK_{j,k}(x)\int_0^t\!\!\int f_k(y,s)\,dy\,ds$.
 \end{remark}
 
\Proof
We introduce a function $\phi(x,t)$ through the relation
\begin{equation}
\label{dcf}
f(x,s)=\biggl(\int f(y,t)\,dy\biggr)g(x)+\phi(x,t),
\end{equation}
where $g_t(x)=(4\pi t)^{-d/2}\exp(-|x|^2/(4t))$ is the gaussian and $g=g_1$.
Identity~\eqref{dcf} is inspired  by a paper by
Schonbek, \cite{Sch91}. See also~\cite{BraV07},
where a similar idea was used to write an asymptotic expansion of the nonlinear term.
From $e^{t\Delta}g=g_{t+1}$ and from the fact that $\P$ commutes with
the heat kernel we get (applying also Lemma~\ref{asG}) the decomposition
$$ \L(f)(x,t)=\KK(x)\int_0^t\!\!\int f(y,s)\,dy\,ds+\mathcal{R}_1(x,t)+\mathcal{R}_2(x,t),$$
where
$$\mathcal{R}_1(x,t)=|x|^{-d}\int_0^t\Psi(x/\sqrt{t+1-s})\int f(y,s)\,dy\,ds$$
and 
$$\mathcal{R}_2(x,t)=\int_0^t\!\!\int\K(x-y,t-s)\phi(y,s)\,dy\,ds.$$
Using that $|\Psi(x)|\le C|x|^{-1}$ (this is a consequence of Lemma~\ref{asG})
we obtain
$$|\mathcal{R}_1(x,t)|\le C|x|^{-d-1}\int_0^t\sqrt{t+1-s}\,\|f(s)\|_1\,ds\le C\sqrt t|x|^{-d-1}.$$
For treating $\mathcal{R}_2$ we use that 
$$\int\phi(y,s)\,dy=0$$
and apply the Taylor formula.
We thus obtain the decomposition
$\mathcal{R}_2=\mathcal{R}_{2,1}+\dots+\mathcal{R}_{2,4}$,
with
$$\mathcal{R}_{2,1}=-\int_0^t\!\!\int_{|y|\le|x|/2}\int_0^1\nabla_x\K(x-\theta y,t-s)\,d\theta
  \cdot y\phi(y,s)\,dy\,ds,$$
  next
\begin{equation*}
\mathcal{R}_{2,2}+\mathcal{R}_{2,3}
   =\biggl(\int_0^t\!\!\int_{|x-y|\le|x|/2}
     +\int_0^t\!\!\int_{|y|\ge|x|/2,\; |x-y|\ge|x|/2}\biggr)
     \K(x-y,t-s)\phi(y,s)\,dy\,ds,
\end{equation*}
and
\begin{equation*}
\mathcal{R}_{2,4}=-\int_0^t \K(x,t-s)\int_{|y|\ge |x|/2} \phi(y,s)\,dy\,ds.
\end{equation*}

The term $\mathcal{R}_{2,1}$ is bounded using estimate~\eqref{denK} and the estimate
(deduced from~\eqref{f3})
$$ \| \,|x|\phi(s)\|_1\le C(1+s)^{-1/2}.$$
This yields $|\mathcal{R}_{2,1}(x,t)|\le C|x|^{-d-1}\sqrt t$.

The other three terms 
can be treated observing that
$$|\phi(y,s)|\le C(1+|y|)^{-d-2}.$$
Then applying Lemma~\ref{lemlog} we get, for $|x|\ge e\sqrt t$,
$|\mathcal{R}_{2,2}(x,t)|\le Ct|x|^{-d-2}\log(|x|/\sqrt t)\le C|x|^{-d-1}\sqrt t$.
For the $\mathcal{R}_{2,3}$ term, we can observe that
the integrand is bounded by $C|y|^{-2d-2}$. Therefore,
$|\mathcal{R}_{2,3}(x,t)|\le Ct|x|^{-d-2}\le C|x|^{-d-1}\sqrt t$.
The term $\mathcal{R}_{2,4}$ can be estimated in the same way.

\endProof

\section{Main results}
\label{sectionmain}

From the previous Lemma we now deduce the following result:
it completes to the case $f\not=0$ the asymptotic profile constructed
in~\cite{BraV07}.

\begin{theorem}
\label{theo1}
 Let $f$, $a$ and $u$ be as in Proposition~\ref{prop1}. We assume that~$f$ satisfies also condition~\eqref{f3}.
Then (the notation is explained in Remark~\ref{rema2})
\begin{equation}
\label{asu}
u(x,t)=e^{t\Delta}a(x)+\KK(x)\int_0^t\!\!\int f(y,s)\,dy\,ds+\mathcal{R}(x,t),
\end{equation}
for some function $\mathcal{R}$ satisfying
\begin{equation}
\label{br}
\forall\, |x|\ge e\sqrt t\colon\quad
|\mathcal{R}(x,t)|\le C|x|^{-d-1}\sqrt{t},
\end{equation}
where $C>0$ is a constant independent on $x$ and $t$.
\end{theorem}

\begin{remark}
This theorem essentially states that if $t>0$ is fixed and $\int_0^t\!\!\int f(y,s)\,dy\,ds\not=0_{\R^d}$,
then 
$$u(x,t)\simeq e^{t\Delta}a(x)+\KK(x)\int_0^t\!\!\int f(y,s)\,dy\,ds \qquad
\hbox{as $|x|\to\infty$}.$$
\end{remark}

\Proof
Owing to our previous Lemma, 
the only thing that remains to do  is to write $u$ through Duhamel formula~\eqref{6a5}
and show that
$$B(u,u)\equiv \int_0^t\!\!e^{(t-s)\Delta}\P\nabla\cdot(u\otimes u)(s)\,ds$$
can be  bounded by $C|x|^{-d-1}\sqrt t$.
In fact, a stronger estimate will be proved.
Note that the convolution kernel $F(x,t)$ of $e^{t\Delta}\P\nabla$ satisfies
to the same decay estimates of $\nabla\K$ (see~\eqref{denK}).
Therefore, after writing 
$$B(u,u)(x,t)=\int_0^t\!\!\int F(x-y,t-s)(u\otimes u)(y,s)\,dy\,ds,$$
 then splitting
the spatial integral into $|y|\le |x|/2$ and $|y|\ge|x|/2$,
and using that $\|F(\cdot,t-s)\|_1=c(t-s)^{-1/2}$, we get
\begin{equation*}
|B(u,u)|(x,t)\le C|x|^{-d-1}\int_0^t\!\!\|u\|_2^2\,ds
\,\,+\,\,(1+|x|)^{-2d} \sqrt t.
\end{equation*}
When $d\ge3$, by estimate~\eqref{ubu} with $\alpha=0$ and $p=2$, we get
$\int_0^t\|u(s)\|_2^2\,ds\le C(1\wedge t)$, which is enough to conclude.
When $d=2$, we have only $\int_0^t\|u(s)\|_2^2\,ds\le  C\log(1+t)$.
(such facts on the $L^2$-norm on~$u$ also follow
from Schonbek's results~\cite{Schon85, Sch91}).
This yields in particular the required bound 
$$|B(u,u)|(x,t)\le C\sqrt t|x|^{-d-1}.$$

\endProof

\begin{remark}
\label{vam}
The asymptotic profile~\eqref{asu} leads us to study
the homogeneous vector fields of the form:
$$\vec m(x)=\KK(x)\vec{c}, \qquad \vec{c}=(c_1,\ldots,c_d).$$
Notice that $\vec m$ has a zero in $\R^d\backslash\{0\}$ \emph{if and only if} $\vec{c}=\vec 0$.
Indeed, we can limit ourselves to the points $\omega\in\S^{d-1}$
and recalling~\eqref{KK} we see that $\vec m$ vanishes at the point $\omega$ if and only if
$$ d\omega_j\omega\cdot\vec{c}=c_j, \qquad\hbox{for $j=1,\ldots,d$}.$$
Multiplying scalarly with $\omega$  we  get $(d-1)\omega\cdot\vec{c}=0$ for all $\omega\in\S^{d-1}$.
As $d\ge2$, we obtain $\vec{c}=0$.
Applying this observation to our situation we get
$$\inf_{x\not=0} |x|^d\biggl|\KK(x)\int_0^t\!\! \int f(y,s)\,dy\,ds\biggr|>0$$
for all $t>0$ such that  $\int_0^t\!\!\int f(y,s)\,dy\,ds\not= 0$.
\end{remark}

A first interesting consequence of Theorem~\ref{theo1} is the following:

\begin{corollary}
\begin{enumerate}

\item
Let $a$, $f$ and $u$ as in Theorem~\ref{theo1}.
We additionally assume that $a=u(0)$ satisfies $a(x)=o(|x|^{-d})$ as $|x|\to\infty$.
Let $t>0$ be fixed, such that $\int_0^t\!\!\int f(y,s)\,dy\,ds\not=0_{\R^d}$.
Then 
for some locally bounded function $R(t)>0$  ($R(t)$ can blow up for  $t\to\infty$ or $t\to0$)
and all $|x|\ge R(t)$, we have 
\begin{equation}
\label{pcc} 
c_t |x|^{-d}\le |u(x,t)|\le c'_t|x|^{-d}
\end{equation}
for some constant $c_t,c'_t>0$ independent on $x$
\item
({\rm Short time behavior of the flow}) 
If in addition $f\in C(\R^+,L^1(\R^d))$ and $f_0=f(\cdot,0)$ is such that
$\int f_0(y)\,dy\not=0_{\R^d}$, then
the behavior of $u$ at the beginning of the evolution
can be described in the following way:
there exists a time $t_0>0$ such that
for all $0<t<t_0$ and all $|x|\ge R(t)$
\begin{equation*}
 ct |x|^{-d}\le |u(x,t)|\le c't|x|^{-d},
\end{equation*}
for some constants $c,c'>0$ independent on $x$ and~$t$.
\end{enumerate}
\end{corollary}

\Proof
Let 
$$\vec m_t(x)\equiv \KK(x)\int_0^t\!\!\int f(y,s)\,dy\,ds.$$

We  apply the asymptotic expansion~\eqref{asu} for~$u$ and study each term in this expression.
By the non-zero mean condition on~$f$ and Remark~\ref{vam},
$3|\vec{m}_t(x)|$ can be bounded from above and from below as in~\eqref{pcc}, for some $0<c_t<c'_t$.
Estimate~$\eqref{br}$ shows that  $|\mathcal{R}(x,t)|\le c_t|x|^{-d}$ provided $|x|\ge R(t)$ and $R(t)>0$ is taken sufficiently large.

On the other hand, the bound $|g_t(x-y)|\le C|x-y|^{-d-1}\sqrt t$ and $\|g_t\|_1=1$ for the heat kernel
lead to
\begin{equation}
\label{ncomp}
\begin{split}
 |e^{t\Delta}a(x)|& \le C\int_{|y|\le |x|/2} \sqrt t |x-y|^{-d-1}|a(y)|\,dy +\int_{|y|\ge |x|/2} g_t(x-y)|a(y)|\,dy\\
&\le C|x|^{-d-1}\sqrt t\|a\|_1+ \hbox{ess\,sup}_{|y|\ge |x|/ 2} |a(y)|
\end{split}
\end{equation}
and since $|a(y)|=o(|y|^{-d})$ as $|y|\to\infty$, this expression is also bounded by $c_t|x|^{-d}$ for $|x|\ge R(t)$ and a sufficiently large $R(t)$.

In conclusion the first and third terms in expansion~\eqref{asu} can be absorbed by $\vec{m}_t(x)$.
The first conclusion of the corollary then follows.
The second conclusion is now immediate, because
for $t>0$ small enough and $|x|\ge R(t)$ large enough,
one can find two constants $\alpha,\beta>0$
such that
$$\alpha t|x|^{-d}\le |\vec m_t(x)|\le \beta t|x|^{-d}$$
(the condition $\int f_0(y)\,dy\not=0_{\R^d}$ is needed only for the lower bound).

\endProof

\medskip

\begin{remark}
The result for the free Navier--Stokes equation ($f\equiv0$)
is different (see \cite{BraV07}): instead of~\eqref{pcc} 
one in general obtains estimates of the form
\begin{equation*}
 c_t |x|^{-d-1}\le |u(x,t)|\le c'_t|x|^{-d-1}.
\end{equation*}
Therefore, an external force, even if small and compactly supported in space-time,
has the effect of increasing the velocity of the fluid particles {\it at all points\/}
at large distances.
\end{remark}

We obtain, as another
consequence of Theorem~\ref{theo1}, the following result
for the large time behavior.

\begin{theorem}
\label{theo2}
Let $a$, $f$ and $u$ as in Theorem~\ref{theo1}.
We additionally assume that $a=u(0)$ satisfies $a(x)=o(|x|^{-d})$ as $|x|\to\infty$
and that 
$$\int_0^\infty\!\!\!\int f(y,s)\,dy\,ds\not=0_{\R^d}.$$
Then, for $\alpha\ge0$, $1<p<\infty$, for some $t_0>0$, some constants $c,c'>0$ and and all $t>t_0$ we have, 
\begin{equation}
\label{wep}
ct^{-\frac{1}{2}(d-\alpha-\frac{d}{p})}\le \Bigl\|(1+|x|)^{\alpha}u(t)\Bigr\|_{p}\le c't^{-\frac{1}{2}(d-\alpha-\frac{d}{p})},
\qquad
\hbox{when $\alpha+d/p<d$}.
\end{equation}
The above inequalities remain true in the limit case $(\alpha,p)=(d,\infty)$.

On the other hand,
\begin{equation}
\label{wep3}
\Bigl\|(1+|x|)^{\alpha}u(t)\Bigr\|_{p} = +\infty \qquad\hbox{when $\alpha+d/p\ge d$},
\end{equation}
and the above equality remains true in the limit cases $p=1$ or ($\alpha>d$ and $p=\infty$).

\end{theorem}

\Proof
The upper bound in~\eqref{wep} has been already proved in Proposition~\ref{prop1}.
The proof makes use of an argument used before in~\cite{BraV07}.
By our assumption on $f$ and Remark~\ref{vam}, the
homogeneous function 
$$\vec m(x)\equiv\KK(x)\int_0^\infty\!\! \int f(y,s)\,dy\,ds$$
does not vanish.
Therefore, for some $c_0>0$,
we obtain $|\vec{m}(x)|\ge c_0|x|^{-d}$.
Let us apply the profile~\eqref{asu}, writing the second term on the right hand side as
$\vec m(x)-\KK(x)\int_t^\infty\!\!\int f(y,s)\,dy\,ds$.
We get for a sufficiently large $M>0$,
all $t>M$ and all $|x|\ge M\sqrt t$, 
 $$|u(x,t)-e^{t\Delta}a|\ge |\vec m(x)|-C|x|^{-d-1}\sqrt t-\textstyle\frac{c_0}{2}|x|^{-d}
    \ge \textstyle\frac{c_0}{3}|x|^{-d}.$$
On the other hand the computation~\eqref{ncomp} guarantees that we can bound, for $|x|\ge M\sqrt t$ and $t>M$,
$|e^{t\Delta}a(x)|\le \frac{c_0}{12}|x|^{-d}$.
We thus get
\begin{equation} 
\label{plww}
|u(x,t)|\ge \textstyle \frac{c_0}{4}|x|^{-d}, \qquad \hbox{for all $|x|\ge M\sqrt t,\quad t>M$}.
\end{equation}
Let $1< p\le\infty$.
Multiplying this inequality by the weight $(1+|x|)^{\alpha p}$,
then integrating  with respect to $x$ on the set $|x|\ge M\sqrt t$, 
we immediately deduce
$$\Bigl\| (1+|x|)^\alpha u(x,t)\Bigr\|_{p} \ge c t^{-\frac{1}{2}(d-\alpha-\frac{d}{p})}$$
for some $c>0$ and all $t>M$.
In the same way, inequality~\eqref{plww} also implies the lower bound~\eqref{wep}
in the limit case $p=\infty$ and  conclusion~\eqref{wep3}.
The upper bounds for~$u$ obtained in Proposition~\ref{prop1} complete the proof.

\endProof

\bigskip

The lower bounds obtained in~\eqref{wep}
are invariant under the time translation $t\rightarrow t+t_*$,
but the hypothesis made on $f$, $\int_0^\infty\!\!\int f(y,s)\,dy\,ds\not=0_{\R^d}$, is not.
The explanation is that the  conditions on the datum,
$a(x)=o(|x|^{-d})$ as $|x|\to\infty$, and   $a\in L^1(\R^d)$, in general, are not conserved at later times.
In fact, according to~\eqref{asu}, one has $|u(x,t_*)|=o(|x|^{-d})$ as $|x|\to\infty$
if and only if $\int_0^{t_*}\!\!\int f(y,s)\,dy\,ds=0_{\R^d}$. 

\medskip
We did not treat the case of external forces with vanishing integrals.
In that case, the principal term in the asymptotic expansion~\eqref{asLf}
of $u$ as $|x|\to\infty$ disappears.
A similar method, however, can be used to write the next
term in the asymptotics  which equals 
$$\nabla\KK(x)\colon\int_0^t\!\!\int y\otimes f(y,s)\,dy\,ds$$
(more explicitly, $\sum_{h,k}\partial_h\KK_{j,k}(x)\int_0^t\!\!\int y_h f_k(y,s)\,dy\,ds$,
for the $j$-component, $j=1,\ldots,d$).
For this we need to put more stringent assumptions on the decay of~$f$:
the spatial decay must be increased of a factor $(1+|x|)^{-1}$
and the time decay by a factor $(1+t)^{-1/2}$. Assumptions~\eqref{f2}-\eqref{f3}
should also be sharpened accordingly.
For well localized initial data $a(x)$, then one would deduce,
bounds of the form
$$ c_t|x|^{-d-1}\le |u(x,t)|\le c'_t|x|^{-d-1}$$
and 
$$\|(1+|x|)^\alpha u(t)\|_p \simeq ct^{-\frac{1}{2}(d+1-\alpha-d/p)},
\qquad \hbox{for $\alpha+d/p<d+1$}.$$
However, as for the free Navier--Stokes equations (see \cite{BraV07}), suitable additional
{\it non-symmetry\/} conditions on the flow and on the matrix $\int_0^t\!\!\int y\otimes f(y,s)\,dy\,ds$
should be added for the validity of the lower bounds.

\section{Acknowledgments}
The preparation of this paper was supported by EGIDE, through the program Huber-Curien ``Star''
N.~16560RK

\bibliographystyle{amsplain}

\end{document}